\documentclass[12pt,a4paper,final]{amsart}
\usepackage[english]{babel}
\usepackage{latexsym}
\usepackage[all]{xy}
\usepackage{amssymb}
\usepackage{amscd}
\usepackage{mathrsfs}
\usepackage[cp850]{inputenc}
\usepackage{url}
\usepackage[notcite,notref]{showkeys}
\usepackage{marginnote}
\theoremstyle{plain}
\newtheorem{nada}{}[section]
\theoremstyle{remark}

\newcommand{\K}{\mathbb{K}}
\newcommand{\Ker}{\operatorname{Ker}}
\newcommand{\cok}{\operatorname{cok}}
\newcommand{\Cok}{\operatorname{Cok}}

\newcommand{\W}{\bf W}
\newcommand{\B}{\ensuremath{\bf B}}
\newcommand{\A}{\ensuremath{\bf A}}

\newcommand{\hM}{\ensuremath{\bf hM}}
\newcommand{\To}{\longrightarrow}

\def\Hom{\operatorname{Hom}}
\def\Ext{\operatorname{Ext}}

\def\PB{\operatorname{PB}}
\def\PO{\operatorname{PO}}

\def\IM{\operatorname{Im}}

\newdir{(>}{{}*!/-10pt/\dir{>}}
\newcommand{\monicar}{\ar@{(>->}}

\def\Ext{\operatorname{Ext}}
\def\Hom{\operatorname{Hom}}
\thanks{Supported in part by MINCIN, Project PID2019-103961GB-C21,  Junta de Extremadura, Projects IB-20038 and IB-18087 and ERDF}

\begin{document}

\title[Quotient Banach spaces after Wegner]{On the category of quotient Banach spaces after Wegner}

\author{F\'elix Cabello S\'anchez}
\author{Jos\'e Navarro Garmendia}
\address{Departamento de Matem\'aticas and IMUEx, Universidad de Extremadura\\
Avenida de Elvas\\ 06071-Badajoz\\ Spain}
\email{fcabello@unex.es\\ navarrogarmendia@unex.es}

\thanks{AMS Subject Class. 2020: 46M18, 18A20}

\noindent{\footnotesize{Version March 4, 2022}}

\bigskip

\bigskip

\maketitle

\begin{abstract}
We study Waelbroeck's category of  Banach quotients after Wegner, focusing on its basic homological and functional analytic properties.
\end{abstract}

\section{Introduction}
Anyone interested in Banach spaces {\em and\,} in category theory knows that {\bf B}, the category of Banach spaces and (linear, bounded) operators, is a nice additive (even quasi-abelian) category that fails to be abelian. The root of this misbehavior of {\bf B} is that, while the categorical kernel of an operator $f:X\To Y$ is  $f^{-1}(0)=\{x\in X:f(x)=0\}$, with the canonical inclusion into $X$, the cokernel of $f$ is given by (the obvious operator from $Y$ onto) $Y/\overline{f(X)}$, where $f(X)=\{f(x):x\in X\}$ and the bar denotes closure in $Y$.

 The inexorable conclusion is that
an operator  $f:X\To Y$ can be monic and epic
 without being an isomorphism in $\B$: just consider any injective operador with dense range but not surjective.

Cain's mark is also responsible for the lack of projective and injective objects in {\bf B}. This statement may be controversial because in Banach space theory we are all taught that 
$\ell_\infty$ is an ``injective Banach space'', that $\ell_1$ is a ``projective Banach space'' and that each finite-dimensional space, equivalently, the ground field, is both.

However the notions between quotes refer to the extension (respectively, lifting) of operators from {\em closed} subspaces, equivalently, through isomorphic embeddings (respectively, quotient maps), while the corresponding categorical notions are much stronger since the corresponding extension (respectively, lifting) must exist through any injective operator (respectively, operator with dense range), and so the only injective or projective object in $\bf B$ is $0$; see  \cite{p} and \cite[Remark 2.1.5]{buhler}.

It is a commonplace, which we share, that homological techniques work fine in abelian categories and that they
can only be adapted to other categories after much work and many sacrifices. (Admittedly, the situation is not that dramatic because we have exact categories: the traditional definition of ``injective\,/\,projective'' objects in Banach space theory correspond to the so-called ``quasi-injective\,/\,quasi-projective'' objects in quasi-abelian categories and to ``$\mathscr E$-injective\,/\,$\mathscr E$-projective'' objects when $\mathscr E$ is the maximal exact structure in $\B$, see \cite{bu-exact, FS, HSW}).

This leads naturally to consider abelianizations of $\B$ --- abelian categories where $\B$ embeds as a full subcategory.

The first such embedding was proposed by Waelbroeck, by using a category of formal quotients, and then by No\"el by means of a category of functors; see \cite[Chapter 2]{w-n} and the references therein, especially \cite{w1, noe}. These constructions were initially motivated by
``practical'' needs arising from holomorphic functional calculus and are more concerned with bornologies than with topologies.

Very recently Wegner \cite{heart} generalized and clarified Waelbroeck ideas, purging them of some errors and flaws; see also \cite[1.3.24]{bbd}, \cite[Chapter IV,\S~2.6]{buhler}, \cite[Part II, \S~14]{HHI}. In this note we  combine Wegner's approach and rather elementary functional analysis to explore some salient features of the Waelbroeck--Wegner category $\W$, mainly of a ``homological'' nature. The paper does not contain any profound results, but confirms that this category is surprisingly (?) well suited for the study of extensions of Banach spaces.

We want to emphasize that, although $\W$ has a long history and pops up in different settings (in \cite{bbd}, as the heart of the natural $t$-structure on the derived category of $\B$, or in \cite{sch}, where Schneiders proves that $\W$ and $\B$ are derived-equivalent), Wegner's ideas are crucial for us because they make $\W$ accesible to explicit computations without the need to introduce a derived category first. We owe this observation (as well as the idea of including the current Section~\ref{sec:char-W}) to the referee of an earlier version of this paper.

Even more recent is Clausen--Scholze's embedding via sheaves, see \cite{scholze}. Curiously enough, Waelbroeck himself  claims to have been inspired by Godement's {\em Topologie alg\'ebrique et th\'eorie des faisceaux}, while No\"el's construction in \cite{noe} follows the Eilenberg--Mac Lane program to the letter.

\section{The Waelbroeck--Wegner category}

In this section we introduce our working category just fixing $\B$ as the ``seed category'' in Wegner's construction. The resulting category corresponds to the Waelbroeck category of quotient Banach spaces (after using the correct definition of  pseudo-isomorphism, see \cite[Definition 2.13.1]{w-n} and compare with \cite[Definition 3]{w1}) and will be denoted by  $\W$. We will be mostly concerned with those properties of $\W$ that depend on the specific features of $\B$ and refer to \cite{heart} for everything else.

\subsection*{Categories of monics}
Let us begin with $\bf M$, the category of monomorphisms in $\bf B$. 
The objects of $\bf M$ are injective operators acting between Banach spaces.
We represent a typical object of $\bf M$ by $f:X'\rightarrowtail X$. A morphism from $f:X'\rightarrowtail X$ to $g:Y'\rightarrowtail Y$ is a
commutative square
\begin{equation}\label{eq:u:f->g}
\xymatrix{
X'  \monicar[d]_f\ar[r]^{u'} & Y' \monicar[d]^g\\
X \ar[r]^{u} & Y\\
}
\end{equation}
in $\B$. The operator $u':X'\To Y'$ is determined by $u$ and so we speak of $u$ as a morphism from $f$ to $g$ and write $u:f\To g$.

The category $\bf hM$ is the homotopical version of $\bf M$: it has the same objects as $\bf M$, with
$$
\Hom_{\bf hM}(f,g)=\frac{\Hom_{\bf M}(f,g)}{J(f,g)},
$$
where $J(f,g)$ is the linear subspace of those morphisms $u:f\To g$ for which there exists an operator $r:X\To Y'$ making the following diagram  commutative
$$
\xymatrix{
X'  \monicar[d]_f \ar[r]^{u'} & Y' \monicar[d]^g\\
X \ar[ur]^r \ar[r]^{u} & Y\\
}
$$
Note that $u=gr\implies u'=rf$ since $g$ is monic.
Clearly, $\hM$ is additive.
The following diagram represents the kernel and cokernel of $u:f\To g$ in ${\bf hM}$: 
\begin{equation}\label{eq:ker cok en hM}
\xymatrixcolsep{3pc}
\xymatrix{(u')^{-1}(0) \monicar[d]_{\text{restriction}} \ar[r] &
X'  \monicar[d]_f \ar[r]^{u'} & Y' \monicar[d]^g \ar[r] & u(X)+Y' \monicar[d]^{\text{inclusion}}\\
u^{-1}(0)  \ar[r] & X  \ar[r]^{u} & Y \ar@{=}[r] &Y\\
}
\end{equation}
where $u(X)+Y'$ carries the norm $z\longmapsto\inf\{\|x\|_X+\|y'\|_{Y'}:z=u(x)+y'\}$ and 
$(u')^{-1}(0)$ and $u^{-1}(0)$ the restrictions of the norms of $X'$ and $X$, respectively \cite[\S~2.2.4]{w-n}. In particular:

\begin{nada}\label{monic-epic} 
$u$ is monic in ${\bf hM}$ if and only if $u^{-1}(0)\subset f(X')$. It is epic if and only if $Y= u(X)+g(Y')$.
\end{nada}

We embed  $\B$ as a full subcategory of $\hM$ sending each Banach space $X$ to the inclusion $0\rightarrowtail X$ (the action on operators is obvious). One has:

\begin{nada}\label{objects_cokers}  Any object $f:X'\rightarrowtail X$ in $\bf hM$ is the cokernel of a morphism between Banach spaces:
\begin{equation*}\label{formalQuotient}
\xymatrix
{0  \monicar[d] \ar[r] & 0 \ar[r] \monicar[d] & X' \monicar[d]^f \\
X'  \ar[r]^{f} & X \ar@{=}[r] & X
}
\end{equation*}
\end{nada}



\begin{nada}
\label{n:w}  $f:X'\rightarrowtail X$ is isomorphic to a Banach space $Y$ (read $0\rightarrowtail Y$) in $\bf hM$ 
 if and only if $f(X')$ is a complemented (hence closed) subspace of $X$ and $X/f(X')$ is linearly homeomorphic to $Y$.
\end{nada}

\subsection*{Exact sequences in a preabelian category}
Exactness will play an increasingly important role as the article progresses, so let us record right now the pertinent definitions. A (finite or infinite) sequence of morphisms
$
\xymatrix{
\cdots\ar[r] & A_{i-1} \ar[r]^-{a_{i-1}} &  A_{i} \ar[r]^-{a_{i}} & A_{i+1} \ar[r] &\cdots
}
$
in an additive category ${\bf A}$ having kernels and cokernels is exact at $A_i$ if ${a_{i}}{a_{i-1}}=0$ and the natural arrow $\Ker \cok {a_{i-1}} \To \Ker a_i$ is an isomorphism:
$$
\xymatrix{
& \Ker(\cok a_{i-1}) \ar[dr] \ar@{..>}[r]&\Ker(a_i)\ar[d]&&\\
\cdots\ar[r] & A_{i-1} \ar[u]^-{\text{induced}} \ar[r]^-{a_{i-1}} &  A_{i} \ar[r]^-{a_{i}} \ar[dr] & A_{i+1} \ar[r] &\cdots\\
&&& \Cok(a_{i-1})\ar[u]_-{\text{induced}}&&
}
$$
The open mapping theorem guarantees that this definition is equivalent to the usual definition for Banach spaces: the kernel of $a_i$ agrees with the range of $a_{i-1}$. A sequence is called exact if it is exact at every position.
A short exact sequence is one of the form
$\xymatrix{\!0\ar[r] & A\ar[r]^\imath & B \ar[r]^\pi & C\ar[r] & 0\!}$.
This means that $\imath$ is a kernel of $\pi$ and $\pi$ is a cokernel of $\imath$; in particular $\imath$ is monic and  $\pi$ is epic.
Sometimes we say that
$\xymatrix{\!A\ar[r]^\imath & B \ar[r]^\pi & C\!}$ is
 {\em short exact} and we say that it splits if $\pi$ admits a right-inverse in ${\bf A}$ or, equivalently, $\imath$ admits a left-inverse in ${\bf A}$. The diagram in {\bf \ref{objects_cokers}} is a short exact sequence in $\hM$.


\subsection*{The category of quotient Banach spaces after Wegner} A morphism $u:f\To g$ in ${\bf hM}$ is called a pulation if (\ref{eq:u:f->g}) is both a pushout and a pullback diagram. As a particularly motivating example, observe that bending a little a short exact sequence of Banach spaces
\begin{equation}\label{eq:YXZ}\xymatrix{
0\ar[r] & Y\ar[r]^\imath & X\ar[r]^\pi & Z\ar[r] & 0
}
\end{equation}
we obtain the pulation 
\begin{equation}\label{associatedPulation}
\xymatrix
{Y \monicar[d]_{\imath} \ar[r] & 0 \monicar[d] \\
 X \ar[r]^{\pi} & Z }
\end{equation}
which is an isomorphism in $\bf hM$ if and only if (\ref{eq:YXZ}) splits, see {\bf \ref{n:w}}.

Being a pulation is really a property of the class of $u$ in ${\bf hM}$ \cite[Lemma 6]{heart}. Pulations form a localizing class (multiplicative system) in ${\bf hM}$ and thus it makes sense to consider the localization of $\bf hM$  with respect to pulations. The resulting category is the Waelbroeck--Wegner category $\bf W$. The objects of $\bf W$ are those of ${\bf hM}$ and there exists a functor $Q:{\bf hM}\To{\bf  W}$ which is the identity on objects and 
takes each pulation into an isomorphism. Every morphism $\phi: f\To g$ in ${\bf  W}$ can be represented as $\phi=Q(u)Q(s)^{-1}$, where $u: h\To g$ is a morphism of ${\bf hM}$ and $s:h\To f$ is a pulation:
\begin{equation*}
\xymatrixrowsep{1pc}
\xymatrix{
&  h \ar[dl]_-s \ar[dr]^-u \\
f && g
}
\end{equation*}
This is called a  ``left roof''. Every morphism can be represented {\em also} by a ``right roof'':
\begin{equation*}
\xymatrixrowsep{1pc}
\xymatrix{
&  k  \\
f \ar[ur]^v   && g \ar[lu]_-t
}
\end{equation*}
See Mili\v{c}i\'c's notes \cite[Chapter 1]{milicic}  for details, including when two of those roofs induce the same morphism in $\bf W$ and \cite[Theorem 10]{heart} for a proof that $\W$ is abelian. (Waelbroeck's proof can be seen in \cite[\S~2.3]{w-n}.)   It is obvious that $\W$ contains $\B$ as a full subcategory.

The following remark implies that the functor $Q:\hM\To \W$ is injective on morphisms. The proof uses quite specific properties of pulations.

\begin{nada}\label{Qu=0}
If $u$ is a morphism in $\hM$ then $u=0$ in $\hM$ if and only if $Q(u)=0$ in $\W$.
\end{nada}

\begin{proof}
One direction is trivial. For the ``if'' part, if $Q(u)=0$ in $\W$ there exists a pulation $s$ such that $su=0$ in $\hM$, see \cite[2.1.4. Lemma]{milicic}. This is depicted in the solid part of the diagram
$$
\xymatrixcolsep{3pc}
\xymatrixrowsep{3pc}
\xymatrix{
X'  \monicar[d]_f \ar[r] & Y' \ar[r] \monicar[d]^>>>>g & Z' 
 \monicar[d]^h\\
X \ar@{..>}[ur]^{r'} \ar[urr]^>>>>>>>>>r \ar[r]^{u} & Y \ar[r]^s & Z\\
}
$$
Since the right square is a pullback diagram the universal property of the pullback applied to the operators $r: X\To Z', u:X\To Y$ yields $r': X\To Y'$ such that $u=gr'$ and shows that the morphism represented by the left square is $0$ in $\hM$. 
\end{proof}

An immediate consequence is:

\begin{nada}\label{cokerQ}
If $u$ is a morphism of $\hM$ then $\ker Q(u)=Q(\ker u)$ and 
$\cok Q(u)=Q(\cok u)$.  
\end{nada}

Thus, applying $Q$ to Diagram (\ref{eq:ker cok en hM}) one obtains explicit descriptions of the kernel and cokernel of $Q(u)$ in $\W$. In particular:

\begin{nada}\label{n:monic-and-epic}
$u$ is monic (respectively, epic) in ${\hM}$ if and only if $ Q(u)$ is monic (respectively, epic) in $\W$.
\end{nada}

A useful criterion to detect if a commutative square formed by arbitrary operators is a pulation (in the sense that is simultaneously a pushout and a pullback) is:

\begin{nada}\label{pula}
A commutative square formed by arbitrary operators in $\bf B$ 
$$
\xymatrix{
A \ar[d]_\gamma \ar[r]^\alpha & B\ar[d]^\beta\\
C \ar[r]^\delta & D
}
$$
is a pulation if and only if the sequence $
\xymatrix{0\ar[r] & A \ar[r]^-{(\alpha,-\gamma)} & B\oplus C \ar[r]^-{\beta\oplus \delta}  \ar[r] & D \ar[r] & 0}$ is exact.
\end{nada}

As one may guess, $(\alpha,-\gamma)(a)=(\alpha(a),-\gamma(a))$ and 
$(\beta\oplus \delta)(b,c)=\beta(b)+\delta(c)$. This can be proved as  \cite[Proposition 2.12]{bu-exact}, taking advantage of the simplifications provided by the peculiarities of $\bf B$.

We can also characterize those pulations that already define isomorphisms in $\hM$:

\begin{nada}\label{pul=isom}
A pulation is an isomorphism in $\hM$ if and only if the associated sequence splits. 
\end{nada}

\begin{proof}
Let $s$ be a pulation from $\imath :X' \rightarrowtail X$ to $\jmath:Y' \rightarrowtail Y$, where $\imath, \jmath$ are assumed to be plain inclusions, so that $s'=s|_{X'}$. The associated sequence is
\begin{equation}\label{assoc}
\xymatrixcolsep{3pc}
\xymatrix{0\ar[r] & X' \ar[r]^-{(\imath,-s)} & X\oplus Y' \ar[r]^-{s\oplus \jmath}  \ar[r] & Y \ar[r] & 0}
\end{equation}
We first remark that, being a pulation, $s$ is an isomorphism in $\hM$ if and only if it admits a left-inverse in $\hM$: indeed, if $\chi:\jmath\To \imath $ satisfies $\chi s={\bf I}_\imath$ in $\hM$, then 
 the identity $s\chi={\bf I}_\jmath$ (in $\hM$) follows from {\bf\ref{Qu=0}}, taking into account that $Q(\chi)Q(s)={\bf I}_\imath$ (in $\W$) and that $Q(s)Q(\chi)={\bf I}_\jmath$ (since $Q(s)$ is an isomorphism in $\W$).
 
 $(\Longleftarrow)\,$ 
Assume that the quotient map of (\ref{assoc}) admits  a left-inverse 
which we write as $(\chi,\eta): Y\To X\oplus Y'$, where $\chi:Y\To X$ and $\eta:Y\To Y'$ are operators. We will show that $\chi$ induces a homomorphism $\imath\To\jmath$ in $\hM$ such that $\chi s$ is (equivalent to) the identity of $\imath$. It is clear that $\chi:Y\To X$ is bounded. To verify that  $\chi:Y'\To X'$ is (well-defined and) bounded, let $P={\bf I}_{ X\oplus Y'}- (\chi,\eta)(s\oplus \jmath)$ be the induced projection on  $X\oplus Y'$ (whose range is $(\imath,-s)(X')$), that is,
\begin{equation*}
P(x,y')= (x,y')- (\chi(s(x)+y'),\eta(s(x)+y')=
(x-(\chi(s(x)+y'),y'- \eta(s(x)+y').
\end{equation*}
In particular
$
P(0,y')=(-\chi(y'), y'-\eta(y')),
$ 
from where it follows that $\chi(y')\in X'$ and that 
$\|\chi: Y'\To X'\|\leq \|(\imath,-s)^{-1}\|\|P\|$.

We now prove that the concatenation
$$
\xymatrix{
X'  \monicar[d]_\imath \ar[r] & Y' \ar[r] \monicar[d]^-\jmath & X' 
 \monicar[d]^-\imath\\
X  \ar[r]^{s} & Y \ar[r]^{\chi} & X\\
}
$$
is 
equivalent to the identity of $\imath$, which amounts to showing that $\chi s- {\bf I}_X$ is bounded from $X$ to $X'$ and is obvious after realizing that
$
P(x,0)= (x-\chi(s(x)), \eta(s(x)))
$.

$(\implies)\,$ The argument is clearly reversible: if $\chi:Y\To X$ is a left-inverse of $s$ in $\hM$ with witness operator $r:X\To X'$, then the map $P:X\oplus Y'\To X'$ defined by $P(x,y')=r(x)+\chi(\jmath(y'))$ is a projection along $(\imath, -s)$.
\end{proof}



\begin{nada}\label{arisesasQ(u)}
Every morphism from $f:L \rightarrowtail \ell_1(I)$ to any object $g$ in $\W$ arises as $Q(u)$, where $u:f\To g$ is a morphism of $\hM$.
\end{nada}

\begin{proof}
The preceding item and the lifting property of $ \ell_1(I)$ imply that every pulation whose codomain is $f:L \rightarrowtail \ell_1(I)$ is an isomorphism in $\hM$.
\end{proof}

\section{A bit of algebra and analysis in {\bf W}}

\subsection*{Distinguished objects}
The following remark is  \cite[Lemma 2.1.22]{w-n}:

\begin{nada}
\label{n:w2}
Each object of $\bf W$ is isomorphic to one of the form $L \rightarrowtail \ell_1(I)$ for some set $I$.
\end{nada}

\begin{proof}
Given $f: X' \rightarrowtail X$ take a surjective operator $\pi:\ell_1(I)\To X$ and form the pullback square (in {\bf B})
$$
\xymatrix{
\PB  \monicar[d] \ar[r] & X' \monicar[d]^f\\
\ell_1(I)\ar[r] & X
}
$$
which is easily seen to be a pulation, hence an isomorphism in $\bf W$.
\end{proof}

The following result characterizes the objects that are isomorphic to Banach spaces in $\W$:

\begin{nada} [Compare to {\bf \ref{n:w}}]\label{n:we}
For a monomorphism $f:X'\rightarrowtail X$ between Banach spaces, the following are equivalent (and imply that  $f$ is isomorphic to $X/f(X')$ in $\W$):
\begin{itemize}
\item[(a)] $f(X')$ is closed in $X$.
\item[(b)] $f$ is isomorphic to a Banach space in $\W$.
\item[(c)] Every morphism $\mathbb K\To f$ admits a left-inverse.
\end{itemize} 
\end{nada}

\begin{proof} (a) $\implies$ (b) is now obvious: $\xymatrix{0 \ar[r] &  X' \ar[r]^f & X  \ar[r] & X/ f(X')  \ar[r] & 0}$ is an exact sequence, and its associated pulation
$$
\xymatrixcolsep{3pc}
\xymatrix{
X'  \monicar[d]_f \ar[r] & 0 \monicar[d]\\
X \ar[r] & X/f(X')\\ 
}
$$ is an isomorphism in $\W$. 
The implication  (b) $\implies$ (c) is clear since (c) is satisfied whenever $f$ is a Banach space, by the Hahn-Banach theorem. (c) $\implies$ (a) Assume $f(X')$ is not closed in $X$ and pick a point $p\in \overline{f(X')}\backslash f(X')$. The linear map $c\in\mathbb K\longmapsto cp\in X$ defines a morphism $\mathbb K\To f$ which does not have a left-inverse in $\W$.
\end{proof}

The preceding result shows that the ground field (hence every nonzero Banach space) fails to be injective in $\W$ and also that every object of $\W$ finds its place in a short exact sequence (see the paragraph after {\bf \ref{n:enough}})
\begin{equation}\label{SEX}
\xymatrixcolsep{3.2pc}
\xymatrix{
X' \monicar[d] \ar@{=}[r] & X' \monicar[d]^f \ar[r] &0 \monicar[d]\\
\overline{f(X')}\ar[r]^-{\text{inclusion}} & X \ar[r]^-{\text{quotient}} & X/\overline{f(X')}
}
\end{equation}
in which the ``quotient'' is a Banach space and the ``subspace'' does not have morphisms to $\K$, apart from zero. We shall prove later ({\bf \ref{noinj}}) that $\W$ does not have injectives. And projectives?

\begin{nada}\label{proj<-->l1}
An object is projective in $\W$  if and only if   it is isomorphic to $\ell_1(I)$ for some set $I$.
\end{nada}

\begin{proof}
$(\Longleftarrow)$ We must prove that every morphism whose domain is $\ell_1(I)$ can be lifted through any epimorphism of $\bf W$. Since $\bf W$ is abelian it has pullbacks and they ``preserve'' epics. Therefore, it suffices to see that every every epimorphism whose codomain is $\ell_1(I)$ admits a right-inverse.

Let $\phi: f\To \ell_1(I)$ be an epimorphism in $\bf W$ that we write as $\phi=Q(u)Q(s)^{-1}$, as a left roof, where $u:g\To \ell_1(I)$ is a morphism in $\bf M$ and $s:g \To f$ is a pulation:
\begin{equation*}
\xymatrixrowsep{0.5pc}
\xymatrix{
& Y' \monicar[dd]_g \ar[dl] \ar[dr]\\
X' \monicar[dd]_f &  &0 \monicar[dd]\\
& Y \ar[dl]_-s \ar[dr]^-u&\\
X &  & \ell_1(I)\\
}
\end{equation*}
Since $Q(s)$ is an isomorphism in $\bf W$ it is clear that $Q(u)$ is epic and (by {\bf \ref{n:monic-and-epic}}) $u$ is epic in $\hM$. This implies that $u:Y\To \ell_1(I)$ is surjective and, by the lifting property of $\ell_1(I)$, there is $v:\ell_1(I)\To Y$ such that $uv$ is the identity of $\ell_1(I)$. Applying $Q$ to the composition
$$
\xymatrix{
0 \monicar[d] \ar[r] & Y' \monicar[d] \ar[r] & X' \monicar[d]\\
\ell_1 \ar[r] & Y\ar[r]^s & X
}
$$
gives a right inverse of $\phi$.

$(\implies)$
Each $f: X'\rightarrowtail X$ in $\W$ comes with the short exact sequence
$$
\xymatrix{
0 \monicar[d] \ar@{=}[r] & 0 \monicar[d] \ar[r] &X' \monicar[d]\\
X'\ar[r]^f & X \ar@{=}[r] & X
}
$$
If this sequence splits then, looking at the left-square and taking into account that $\bf B$ is a full subcategory of $\bf W$, we see that $f$ is a topological embedding with complemented range and, by {\bf \ref{n:we}}, $f$ is isomorphic to the Banach space $X/f(X')$ which in its turn is isomorphic to a complemented subspace of $X$. By {\bf \ref{n:w2}} one can assume $X=\ell_1(I)$ and the result follows from K\"othe's remark that  every complemented subspace of $\ell_1(I)$ is isomorphic to  $\ell_1(J)$ for some $J$, see \cite{koethe} --- or \cite{ort} for a proof written in English.
\end{proof}

\begin{nada}\label{n:enough}
$\bf W$ has enough projectives.
\end{nada}

\begin{proof}
This means that every object is the codomain of an epimorphism whose domain is projective, which is clear for objects of the form $L\rightarrowtail \ell_1(I)$ since $\ell_1(I)$ is projective in $\bf W$, hence for all objects, see {\bf \ref{n:w2}}.
\end{proof}

Now, a standard procedure gives rise to resolutions: we may take a surjective operator  $\pi_1:\ell_1 (I_1) \longrightarrow L$; then take a surjective operator  $\pi_2:\ell_1 (I_2) \longrightarrow \ker\pi_1 $, and so on. 

\begin{nada}\label{projective_Resolution}
Any object $f: X'\rightarrowtail X $ 
in $\W$ admits a projective resolution
$$
\xymatrix{
\ldots \ar[r] & 0 \monicar[d] \ar[r] & 0 \monicar[d] \ar[r] & 0 \monicar[d] \ar[r] & 0 \monicar[d] \ar[r] & X' \monicar[d]^f \\
\ldots \ar[r]  & \ell_1 (I_3) \ar[r]^-{\pi_3} & \ell_1(I_2) \ar[r]^-{\pi_2} & \ell_1 (I_1) \ar[r]^-{\pi_1} & \ell_1 (I_0) \ar[r]^-{\pi_0} & X
}
$$ where the maps ${\pi_i}$ have closed range for $i\neq 1$, the range of $\pi_1$ equals $L=\pi_0^{-1}(f(X'))$, and $f$ is isomorphic to $L\rightarrowtail \ell_1 (I_0)$ in $\W$.
\end{nada}

\subsection*{Credits} Most of the results presented so far are due to Waelbroeck and / or  No\"el (or their direct translations to Wegner's language / the language of Banach spaces).
In any case, there is the small problem (which we will skip) of clarifying the relationship between No\"el's category (which is commonly accepted as being equivalent to Waelbroeck's category of bornological spaces) and ``our'' $\W$.

Diagram (\ref{eq:ker cok en hM}) summarizes the information in \cite[\S2.2.4]{w-n} and {\bf \ref{monic-epic}} corresponds to \cite[Corollaries 2.2.17 and 2.2.19]{w-n}.
{\bf \ref{n:w}} is mentioned in passing in \cite{w1}, see the paragraph closing \S~1.
Assertion {\bf \ref{pul=isom}} generalizes as well as explains  \cite[Lemma 2.1.4]{w-n}. {\bf \ref{arisesasQ(u)}} corresponds to the statement ``Every morphism defined on a standard Q-space is strict'' in \cite[p. 84]{w-n}.

Assertion {\bf \ref{n:w2}} corresponds to \cite[Lemma 2.1.22]{w-n}, while {\bf \ref{n:enough}} is mentioned in the introduction to Chapter 2 of \cite{w-n}. Both {\bf \ref{n:enough}} and the ``if part'' of {\bf \ref{proj<-->l1}} can be attributed to No\"el: they correspond to Propositions 5.6 and 5.2 in \cite{noel-thesis}. However No\"el works in a larger category {\bf QESPC} which must be very different from $\W$. Indeed \cite[Th\'eor\`eme 8.2]{noel-thesis} claims that {\bf QESPC} has enough injectives, which is puzzling since $\W$ presents the opposite behavior:

\begin{nada}\label{noinj}
$\bf W$ has no injective object apart from zero.
\end{nada}

\begin{proof}
The proof of {\bf \ref{n:we}} shows that if the object $f:X'\rightarrowtail X$ is injective (and non-zero) in $\W$, then $f(X')$ must be a proper dense subset of $X$. Otherwise $\K$ is a direct summand of $f$ and we have seen that $\K$ fails to be injective. 
Assume
that $f$ has dense range and is not surjective and let us exhibit a monic whose domain is $f$ that does not admit a left-inverse. Let us form the following pushout diagram in $\B$:
$$
\xymatrix{
X'\ar[r]^f \monicar[d]_f & X \ar[d]^F \\
X\ar[r]^{G} & \PO
}
$$
The operator $F$ is injective (hence an object of $\W$; this is not automatic for pushout squares, see the example in \cite[p. 338]{HSW}). Indeed the pushout space can be constructed as $\PO=(X\oplus X)/\overline{(f,-f)(X')}$ and then $F(x)=(0,x)+\overline{(f,-f)(X')}$ and $G(x)=(x,0)+\overline{(f,-f)(X')}$  which cannot be zero unless $x=0$. Therefore, the square above can be regarded as a monomorphism $G:f\To F$ in $\W$. Assume that $G$ (that is, $Q(G)$) admits a left-inverse $\phi$ in $\W$. Let us represent $\phi$ as a right roof: $\phi= Q(s)^{-1}Q(v)$, where $v$ is a morphism of $\hM$ and $s$ is a pulation. One has ${\bf I}_f=  Q(s)^{-1}Q(v) Q(G)$, thus $ Q(s)=Q(v) Q(G)$ and consequently $s=vG$ (in $\hM$) which means that $vG$ is a pulation. If $g:Y'\rightarrowtail Y$ is the common codomain of $v$ and $s$, then the outer rectangle of the commutative diagram 
$$
\xymatrix{
X'\ar[r]^f \monicar[d]_f & X \monicar[d]^F \ar[r]^{v'}& Y' \monicar[d]^g \\
X\ar[r]^{G} & \PO \ar[r]^v & Y
}
$$
is a pulation and in particular (see {\bf \ref{pula}}) the operator
$(f, -v'f): X'\To X\oplus Y'$ must have closed range. Pick $x\in X\backslash f(X')$ and note that $(x, -v'(x))$ cannot be in the range of $(f, -v'f)$. If $(x_n')$ is a sequence in $X'$ such that $f(x_n')\To x$ in $X$ it is clear that $(f, -v'f)(x_n')$ converges to  $(x, -v'(x))$ and we have reached a contradiction.
\end{proof}

\subsection*{Extensions}
Let ${\bf A}$ be an additive category having kernels and cokernels and let $A,B$ be objects of $\bf A$. An $n$-extension of $A$ by $B$ is an exact sequence
\begin{equation}\tag{$\mathscr A$}
\xymatrix{
0 \ar[r] & B \ar[r] &  A_{1} \ar[r] & A_{2} \ar[r] &\cdots \ar[r] & A_n \ar[r] & A \ar[r] & 0
}
\end{equation}
Within the class of $n$-extensions of $A$ by $B$, one considers the least relation of equivalence that makes $(\mathscr A)$ equivalent to 
\begin{equation*}
\xymatrix{
0 \ar[r] & B \ar[r] &  B_{1} \ar[r] & B_{2} \ar[r] &\cdots \ar[r] & B_n \ar[r] & A \ar[r] & 0
}
\end{equation*}
if there is a commutative diagram
$$
\xymatrix{
0 \ar[r] & B \ar@{=}[d] \ar[r] &  A_{1} \ar[d] \ar[r] & A_{2} \ar[d] \ar[r] &\cdots \ar[r] & A_n \ar[d] \ar[r] & A \ar@{=}[d] \ar[r] & 0\\
0 \ar[r] & B \ar[r] &  B_{1} \ar[r] & B_{2} \ar[r] &\cdots \ar[r] & B_n \ar[r] & A \ar[r] & 0
}
$$
Let $\Ext^n_{\bf A}(A,B)$ denote the {\em class} of those equivalence classes. If $\bf A$ is either $\bf B$ or $\bf W$ then $\Ext^n_{\bf A}(A,B)$ is a set that can be given a natural linear structure by using pullbacks and pushouts. See \cite[Chapter~6]{FS} for a good exposition in the setting of {\em exact categories} --- note that these include both $\W$ (which is abelian) and $\B$ (since it is quasi-abelian). The interested reader can find the translation of some classical results on Banach spaces to the language of $\Ext^n_{\bf B}$ in \cite{CCG}.


We can describe $n$-extensions in $\bf W$ by using projectives in the standard way: pick objects $a:A'\rightarrowtail A$ and $b:B' \rightarrowtail B$ and use {\bf \ref{n:enough}} and the abelian character of $\W$ to get an exact sequence
\begin{equation}\tag{$\mathscr P$}
\xymatrixcolsep{2.5pc}
\xymatrix{
0 \ar[r] & k \ar[r]^-\varkappa &  p_{n} \ar[r]^-{\pi_n} & p_{n-1} \ar[r]^-{\pi_{n-1}} &\cdots \ar[r] & p_1 \ar[r]^{\pi_1} & a \ar[r] & 0
}
\end{equation}
with each $p_i$ projective in $\bf W$ (that is, isomorphic to $0\rightarrowtail\ell_1(I_i)$ for some set $I_i$). If $\phi:k\To b$ is a morphism in $\bf W$ then, taking pushouts successively, one obtains the commutative diagram
\begin{equation}\label{eq:POsuc}
\xymatrix{
0 \ar[r] & k \ar[r]^-\varkappa  \ar[d]_\phi&  p_{n} \ar[d]\ar[r] & p_{n-1} \ar[d]\ar[r] &\cdots \ar[r] & p_1 \ar[d] \ar[r] & a\ar@{=}[d] \ar[r] & 0\\
0 \ar[r] & b \ar[r] &  \PO_{1} \ar[r] & \PO_{2} \ar[r] &\cdots \ar[r] & \PO_{n} \ar[r] & a \ar[r] & 0
}
\end{equation}
in which the lower sequence is exact --- and belongs to the class of zero if and only if $\phi$ ``extends'' to $p_n$, with the meaning that there exists a morphism $\Phi: p_n\To b$ such that $\phi=\Phi\varkappa$, although this point is not completely trivial. Actually every $n$-extension
$$
\xymatrixcolsep{2.5pc}
\xymatrix{
0 \ar[r] & b \ar[r]^-\beta &  a_{n} \ar[r]^{\alpha_n} & a_{n-1} \ar[r]^{\alpha_{n-1}} &\cdots \ar[r] & a_1 \ar[r]^{\alpha_1} & a \ar[r] & 0
}
$$
arises in this way, up to equivalence,
as can be seen looking carefully at the following diagram
$$
\xymatrixcolsep{1.5pc}
\xymatrixrowsep{0.75pc}
\xymatrix{
&&&\ker\pi_{n-1} \ar[dd] \ar[rd] &&  \ker\pi_{2} \ar[dd] \ar[rd] &&  \ker\pi_{1} \ar[dd] \ar[rd] \\
0 \ar[r] & k \ar[dd] \ar[r]^-\varkappa &  p_{n} \ar[dd] \ar[ru]  \ar[rr] & & \cdots  \ar[ru]  \ar[rr] & & p_2 \ar[dd] \ar[ru]  \ar[rr] & & p_1 \ar[dd] \ar[r] & a \ar@{=}[dd] \ar[r] & 0\\
&&&\ker\alpha_{n-1} \ar[rd] &&  \ker\alpha_{2} \ar[rd] &&  \ker\alpha_{1} \ar[rd] \\
0 \ar[r] & b \ar[r]^-\beta &  a_{n} \ar[ru]  \ar[rr]^{\alpha_n} & & \cdots  \ar[ru]  \ar[rr] & & a_2  \ar[ru]  \ar[rr]^{\alpha_2} & & a_1 \ar[r]^{\alpha_1} & a \ar[r] & 0
}
$$
and using the universal property of the pushout. All of this is usually summed up by the familiar mantra ``Yoneda $\Ext^n$ is naturally equivalent to $\Ext^n$ via projectives'', see \cite[Chapter 6]{FS} or \cite[Chapter IV]{HS}.

The operators in any exact sequence of Banach spaces beginning and ending by zero have closed range and {\bf \ref{n:cokB=cokW}} implies that every $n$-extension of Banach spaces can be regarded also as an $n$-extension in $\bf W$. Thus, given Banach spaces $A, B$ and $n\geq 1$, we have a mapping $\Ext^n_{\bf B}(A,B)\To \Ext^n_{\bf W}(A,B)$ which becomes (is) a natural transformation of functors in the appropriate setting. It is in fact a natural equivalence:

\begin{nada}
If $A,B$ are Banach spaces then $\Ext_{\bf B}^n(A,B)=\Ext_{\bf W}^n(A,B)$ 
for all $n\geq 1$.
\end{nada}

\begin{proof}
We need to prove two things: that if $A,B$ are Banach spaces every $n$-extension of $A$ by $B$  in $\bf W$ is equivalent to one coming from $\bf B$ (which means that $\Ext^n_{\bf B}(A,B)\To \Ext^n_{\bf W}(A,B)$ is surjective) and that if two $n$-extensions of Banach spaces become equivalent in $\Ext^n_{\bf W}(A,B)$ then they were equivalent in $\Ext^n_{\bf B}(A,B)$ (which means that the map above is injective).

Both follow from the preceding observations. If $a$ denotes the ``image'' of $A$  in $\bf W$ (that is, $a$ is $0\rightarrowtail A$) then the presentation ($\mathscr P$) can be taken in $\bf B$ (and is exact there). If $b$ denotes $0\rightarrowtail B$, then the whole Diagram (\ref{eq:POsuc}) lives in $\B$ (since the arrows in $(\mathscr P)$ have closed ranges) and $\Ext^n_{\bf B}(A,B)\To \Ext^n_{\bf W}(a,b)$ is onto. Moreover, as $\B$ is a full subcategory of $\W$, the morphism $\phi$ extends to $p_n$ in $\W$  if and only if it extends in $\B$, which shows that  $\Ext^n_{\bf B}(A,B)\To \Ext^n_{\bf W}(a,b)$ is injective.
\end{proof}

\section{A characterization of $\W$}\label{sec:char-W}

This section explains the singular role that $\W$ plays among the {abelianizations} of $\B$. The main result {\bf \ref{propiedad_universal}} is contained in B\"uhler's \cite[Theorem 2.2.3 (ii)]{buhler}  --- note that ``our'' preceding item is contained in the first part of the same theorem! The proofs (in the particular case of $\B$) are much simpler, partly because of Wegner's direct approach and partly because we can take advantage of the resolutions of {\bf \ref{projective_Resolution}}. 
Before proceeding let us cheat a bit about how $\W$ is constructed from $\B$.

\subsection*{You could have discovered $\W$} 
As we have already mentioned, $\B$ is quasi-abelian, but not abelian, and so  
there are two types of monomorphisms in $\B$: the ``good'' ones (monic arrows that are kernels in $\B$), and the ``bad'' ones (the others).  

This naturally leads to consider abelianizations that amend these bad monomorphisms and, at the same time, preserve the good.
We can make precise these requirements, imposing that the abelianization functor $a:\B\To {\bf A}$:  
\begin{enumerate}
\item[($\flat$)] Preserves monomorphisms: if $f \colon X' \To X$ is  monic, then $af \colon aX' \To aX$ has to be monic, hence a kernel, because ${\bf A}$ is abelian. 

\item[($\sharp$)] Is exact: it preserves short exact sequences.
\end{enumerate}
The condition of fixing the bad monics
is already accomplished by ${\bf hM}$: each monic arrow $f \colon X' \rightarrowtail X$ between Banach spaces is now a kernel,  see {\bf \ref{objects_cokers}}. 
However, ${\bf hM}$ spoils many exact sequences of Banach spaces: if
\begin{equation}\label{eq:spoils}
\xymatrix{
0\ar[r] &
 Y\ar[r]^\imath & X\ar[r]^\pi & Z
 \ar[r] & 0
}
\end{equation}
is an exact sequence in {\bf B}, then 
\begin{equation}\label{exampleFails}
\xymatrix
{0  \monicar[d] \ar[r] & 0 \monicar[d] \ar[r] & 0 \monicar[d] \\
Y  \ar[r]^{\imath} & X \ar[r]^{\pi} & Z
} 
\end{equation}
is short exact in $\hM$ if and only if the square on the right is the cokernel of the square on the left. But the cokernel of the square on the left in ${\bf hM}$ {is} $\imath \colon Y \rightarrowtail X$ {\em itself} and therefore (\ref{exampleFails}) is exact in $\hM$ if and only if the arrow defined by the square
\begin{equation}\label{origenPulaciones}
\xymatrix
{Y \monicar[d]_{\imath} \ar[r] & 0 \monicar[d] \\
 X \ar[r]^{\pi} & Z }
 \end{equation}
  (which is always monic and epic, see {\bf \ref{monic-epic}})
 is an isomorphism in ${\bf hM}$, something which is not unless (\ref{eq:spoils})
 splits (so {\bf \ref{pul=isom}} says). 

To solve this issue, we can {\it force} these maps to be isomorphisms, localizing the category. 
But the smallest multiplicative system on ${\bf hM}$ containing the squares of the type (\ref{origenPulaciones}) is precisely the class of all pulations and so we arrive at $\W$.

\subsection*{The universal property of $\W$.}
Let $w \colon \B \To \W$ be the abelianization functor.
En route to our characterization of $\W$, let us now check that $w$ satisfies conditions 
($\flat$) and ($\sharp$) above:

\begin{nada}\label{n:cokB=cokW}
The embedding of $\B$ into $\W$ is an exact functor that preserves monics.
\end{nada}

\begin{proof}
The part concerning monics is clear since both $\B\To \hM$ and the localization functor $Q:\hM\To\W$ preserve monics. The part concerning exactness is implicit in the preceding discussion: if (\ref{eq:spoils}) is a short exact sequence in $\B$, then (\ref{exampleFails}) is short exact in $\W$ because
\[
\xymatrix
{0  \monicar[d] \ar[r] & 0 \monicar[d] \ar[r] & Y \monicar[d]^\imath \\
Y  \ar[r]^{\imath} & X \ar@{=}[r] & X
} 
\]
 is 
  and (\ref{origenPulaciones}) is an isomorphism in $\W$.
\end{proof}

\begin{nada}\label{propiedad_universal}
Let ${\bf A}$ be an abelian category and let $a\colon \B \To \A $ be an exact functor that preserves monics. 
Then there exists a unique, up to natural equivalence, exact functor $\bar{a} \colon \W \To \A$ such that $ a = \bar{a} \circ w$, as illustrated in the diagram
\begin{equation*}\label{TrianglePU}
\xymatrixrowsep{1pc}
\xymatrix
{\B  \ar[dr]_w \ar[rr]^a &  & \A  \\
 & \W \ar@{..>}[ur]_{\bar{a}} & 
}
\end{equation*}
\end{nada}

\begin{proof} We have already seen ({\bf \ref{objects_cokers}}) that any object $f \colon X' \rightarrowtail X$ in $\W$ can be regarded as the cokernel of $w(f)$. Therefore, if such an exact functor $\bar{a}$ exists, then
$$ \bar{a} \left( f \colon X' \rightarrowtail X \right) \simeq \bar{a} \left( \Cok(wf) \right) \simeq \Cok ( \bar{a} w f) = \Cok (af) \ . $$

Therefore, we define a functor $\bar{a}\colon {\bf M} \To \A$ mapping each object $f \colon X' \rightarrowtail X$ to
a {\em fixed} cokernel $Z_{af}$ of $af$ in ${\bf A}$ (which is assumed to be $aX$ if $f$ is $0 \rightarrowtail X$), and each morphism $u$ in ${\bf M}$
\begin{equation}\label{eq:looking}
\xymatrix{
X'  \monicar[d]_f\ar[r]^{u'} & Y' \monicar[d]^g\\
X \ar[r]^{u} & Y\\
}
\end{equation} to the morphism $\bar{a}u$ in $\A$ given by the dotted arrow in the commutative diagram
\begin{equation}\label{a_barra_de_flecha}
\xymatrix{
aX'  \monicar[d]_{af} \ar[r]  & aY' \monicar[d]^{ag}\\
aX \ar[d] \ar[r]^{au} & aY \ar[d] \\
Z_{af} \ar@{..>}[r]^{\bar{a}u} & Z_{ag} \\
} .
\end{equation}
This functor $\bar{a}$ factors through ${\bf hM}$ because any morphism $u:f\To g$ in $J(f,g)$ is transformed into the null map: if we have commutative diagram of operators in $\B$ 
$$
\xymatrix{
X'  \monicar[d]_f \ar[r]^{u'} & Y' \monicar[d]^g\\
X \ar[ur]^r \ar[r]^{u} & Y\\
}
$$ the corresponding map $\bar{a}u$, fitting in the commutative diagram below, must be the zero map in $\A$ since $z$ is epic:
\begin{equation*}
\xymatrix{
aX'  \monicar[d]_{af} \ar[r] & aY' \monicar[d]^{ag}\\
aX \ar[d]_z \ar[ur]^{ar} \ar[r]^{au} & aY \ar[d] \\
Z_{af} \ar@{..>}[r]^{\bar{a}u} & Z_{ag} \\
} .
\end{equation*}

Also, $\bar{a}$ factors through the localization $Q:\hM\To\W$, because it transforms pulations into isomorphisms. To see this, assume that (\ref{eq:looking}) is a pulation $u:f\To g$ in ${\bf hM}$. 
Looking at this through $a$ in ${\bf A}$ we obtain another commutative diagram
\begin{equation*}
\xymatrix{
aX'  \monicar[d]_{af}\ar[r]^{au'} & aY' \monicar[d]^{ag}\\
aX \ar[r]^{au} & aY
}
\end{equation*}
In view of {\bf \ref{pula}}, 
 the sequence
$
\xymatrixcolsep{2.5pc}
\xymatrix{0\ar[r] & aX' \ar[r]^-{(au'\!,-af)} & aY'\oplus aX \ar[r]^-{ag\oplus au}  \ar[r] & Y \ar[r] & 0}
$ is exact in {\bf A}. By B\"uhler's \cite[Proposition 2.12, (ii)$\implies$(iv)]{bu-exact} the preceding square can be {\em completed} to a commutative diagram
\begin{equation*}
\xymatrix{
aX'  \monicar[d]_{af}\ar[r]^{au'} & aY' \monicar[d]^{ag}\\
aX\ar[d] \ar[r]^{au} & aY \ar[d]\\
C \ar@{=}[r] & C } .
\end{equation*}
with short exact  vertical sequences.
A quick glance at (\ref{a_barra_de_flecha}) reveals that $\bar{a}(u)$ is an isomorphism in ${\bf A}$. 

On the other hand, it is clear that any exact functor $\W\To {\bf A}$ whose ``restriction'' to $\B$ agrees with $a$ is naturally equivalent to $\bar{a}$. 

\smallskip

It only remains to see that $\bar{a}:\W\To {\bf A}$ is exact. To this end we need to compute the (left) derived functors of $\bar{a}$, and prove that $L^0\bar{a}$ is naturally equivalent to $ \bar{a}$ it self (which is not automatic since we do not know {\em a priori} that $\bar{a}$ is right-exact) and that $L^1 \bar{a} = 0$, see \cite[Chapter V, \S 5]{CE}.
 As a consequence, it will follow that $\bar{a}$ is exact: indeed, if $0\To f\To g\To h\To 0$ is exact in $\W$ there is a commutative diagram in ${\bf A}$
 \[
 \xymatrix{
 0  \ar@{=}[d]\ar[r] & \bar{a}(f) \ar[r] \ar[d] & \bar{a}(g) \ar[r]  \ar[d] & \bar{a}(h)  \ar[d] \ar[r] & 0\\
 L^1\bar{a}(f) \ar[r] & L^0\bar{a}(f) \ar[r] & L^0\bar{a}(g) \ar[r] & L^0\bar{a}(h) \ar[r] & 0
 }
 \]
in which the second row is exact (since $L^\bullet\bar{a}$ are derived functors) and the vertical arrows are the isomorphisms given by the natural equivalence between $\bar{a}$ and $L^0\bar{a}$.

Recall that the homology of the complex
\begin{equation}\label{eq:A}\tag{${\mathscr C}$}
\xymatrix{
\cdots \ar[r] & C_{n+1} \ar[r]^{c_{n+1}} & C_n \ar[r]^-{c_n} & C_{n-1} \ar[r] & \cdots
}
\end{equation}
in an abelian category (at the $n$-th spot) is defined as
$
H_n({\mathscr C})=\Cok\big(\operatorname{Im} c_{n+1}\To \Ker c_n \big)
$, where the image of the arrow $c:C\To D$ is $\operatorname{Im}  (c)=\Ker(\cok c)$. Note the factorization
\[
\xymatrixrowsep{0.5pc}
\xymatrix{C \ar[dr]_-{\text{epic}} \ar[rr]^c & & D\\
& \operatorname{Im}  (c) \ar[ur]_-{\text{monic}} }
\]
The value of the $n$-th derived functor $L^n \bar{a}$ on any object $f:X' \rightarrowtail X$   
of $\W$ is 
obtained picking any
 projective resolution of $f$:
\begin{equation*}
\xymatrixcolsep{2.5pc}
\xymatrix{
\cdots \ar[r]&  p_{2} \ar[r]^-{\pi_2} & p_{1} \ar[r]^-{\pi_{1}}  \ar[r] & p_0 \ar[r]^{\pi_0} & f \ar[r] & 0
}
\end{equation*}
and taking 
the homology of the complex 
$
\xymatrixcolsep{1.75pc}
\xymatrix{
\cdots \ar[r]&  \bar{a} p_{2} \ar[r] & \bar{a} p_{1} \ar[r]  \ar[r] & \bar{a} p_0 \ar[r] & 0
}
$ 
at the $n$-th spot.
Of course this can be done using whatever resolutions we like best, so we can use those of {\bf \ref{projective_Resolution}} to conclude that the derived functors of $\bar{a}$ at $f$ are the homology of the complex 
\[
\xymatrix{
\cdots \ar[r]&  a\,\ell_1(I_2)   \ar[r]^{a \pi_2} & a\,\ell_1(I_1)   \ar[r]^{a \pi_1} & a\,\ell_1(I_0) \ar[r] & 0
}
\]
Hence,
$$
L^0 \bar{a} (f)= \Cok\big(\IM(a\pi_1) \To a\,\ell_1(I_0)\big) \simeq 
 \Cok\big( af: aX'\To aX\big)= \bar{a} (f)
$$
on objects and the fact that $L^0 \bar{a} $ is naturally equivalent to $\bar{a} $ follows from the Comparison Theorem \cite[2.2.6, p. 35]{weibel}, while $
L^n \bar{a} (f)=0$ for $n\ge 1$ because the sequence
$$
\xymatrix{
\cdots \ar[r]&  \ell_1(I_2)   \ar[r] & \ell_1(I_1)   \ar[r] & \ell_1(I_0)
}
$$
was exact in $\B$.
\end{proof}

The proof shows that for any abelian category $\A$, composition with $w$ is part of an equivalence between the category of exact functors $\W \To \A$ and the category of exact functors $\B \To \A$ that preserve monics. This immediately implies that the pair ($\W$, $w$) is unique,  up to unique equivalence, among all pairs ($\A$, $a$) of an abelian category and an abelianization satisfying the universal property in {\bf \ref{propiedad_universal}}.

\subsection*{Condensed stuff} An immediate consequence of {\bf \ref{propiedad_universal}} is that the natural embedding of $\B$ into Clausen--Scholze category ${\bf CS}$ induces an exact functor $\W\To {\bf CS}$. Without entering into details, the relevant feature is that the image of a Banach space $X$ in ${\bf CS}$ is the sheaf of continuous functions $\big(C(K,X))_{K}$, where $K$ runs over all zero-dimensional compacta whose topological weight is bounded by an appropriately huge cardinal, and the action on operators is given by composition. As one may guess, the abelianization  $\B\To {\bf CS}$ preserves monics and is exact since, whenever 
$
\xymatrix{
0\ar[r] & Y\ar[r]^\imath & X\ar[r]^\pi & Z\ar[r] & 0
}
$ is an exact sequence of Banach spaces, the induced sequence
$$
\xymatrix{
0\ar[r] &  C(K,Y) \ar[r]^{\imath_*}  & C(K,X) \ar[r]^{\pi_*} & C(K,Z)\ar[r] & 0
}
$$ is again exact. The only doubtful point is whether $\pi_*$ is surjective and this follows from the existence of a continuous section of $\pi$ --- the celebrated Bartle--Graves {\em selection} theorem, see \cite[Proposition 2.19 (ii)]{BL} for a simplified proof.

All this raises the question of whether $\W$ is a full subcategory of ${\bf CS}$ or not.

\section{Ramblings}

\subsection*{A remarkable functor} Let $\bf sB$ be the category of semi-Banach spaces and bounded operators (see \cite[p. xi]{buhler} or \cite[Note 1.8.1]{hmbst}).
There is a functor $S:\W\To {\bf sB}$ sending $f:X'\rightarrowtail X$ to the space $X/f(X')$ seminormed by 
$$
|x+f(X')|_{X/f(X')}= \inf_{x'\in X'}\|x+f(x')\|.
$$
$S$ is the topologized version of the functor $R$ appearing in \cite[Section~4]{heart} and is described in more abstract terms in \cite[Lemma on p.~xi]{buhler}. To see the action of $S$ on morphisms take $g:Y'\rightarrowtail Y$ and $u:f\To g$ in $\bf M$ and define $S(u): X/f(X')\To Y/g(Y')$ by $x+f(X')\longmapsto u(x)+g(Y')$.  It is really easy to see that $S(u)$ is bounded by $\|u\|$ and depends only on the class 
of $u$ in $\bf hM$. To ensure that $S$ acts from $\W$ to $\bf sB$ one only has to see that $S$ carries pulations into isomorphisms, taking then avantage of the universal property of the localization \cite[Theorems 1.1.1 and 2.1.2]{milicic}. But a pulation square in $\B$ 
is also a pushout in $\bf V$, the category of vector spaces and linear maps, and so in $\bf sB$  (which is not true, as a rule, for mere pushout squares):

\begin{nada}
Let 
$$
\xymatrix{
X' \ar[r]^v \ar[d]_f & Y' \ar[d]^g\\
X \ar[r]^u & Y
}
$$
be a pulation square in $\B$, in which no operator is assumed to be injective. If $Z$ is a seminormed space and $v':Y'\To Z, f':X\To Z$ are operators such that $vv'=f'f$ then there is a unique operator $h:Y\To Z$ such that $f'=hu, v'=gh$.
\end{nada}
(Assemble your own proof: the key point is that the set $\{(v(x'), f(x'): x'\in X'\}$ is closed in $Y'\times X$. See why?) Hence a pulation square induces two isomorphisms of seminormed spaces: namely $\overline{u}: X/f(X')\To Y/g(Y')$ and $\overline{g}: Y'/v(X')\To Y/u(X)$.

Our attempts to gain a more complete understanding of $\W$ come up against the hard truth that we are basically unable to distinguish two objects of $\W$ if they have isomorphic  images in $\bf sB$, see {\bf \ref{noisom}} and {\bf \ref{lr}} for the only exceptions at hand.

Note that $S(f)=\big{(}X/\overline{f(X')}\big{)}\oplus \big{(} \overline{f(X')}/f(X')\big{)}$ where the first summand is a Banach space and the second carries the trivial seminorm. Incidentally, the space $\overline{f(X')}/f(X')$ {\em is} the difficult part since it does not allow for much nuance:

\begin{nada}\label{uncountable}
Let $f:X' \rightarrowtail X$ be a monomorphism between Banach spaces. If $\overline{f(X')}/f(X')$ has a (finite or infinite) countable Hamel basis, then $f$ has closed range, so $\overline{f(X')}/f(X')=0$, and $f$ is isomorphic to a Banach space in $\W$.
\end{nada}

\begin{proof}
The proof uses a version the closed graph theorem due to De Wilde \cite{dewilde}.

Let $(h_i)_{i\in I}$ be a Hamel basis of $\overline{f(X')}/f(X')$ with $I$ countable and, for each $i$, pick $e_i\in X$ such  that $h_i=e_i+f(X')$. Let $\K^{(I)}$ be the (topological linear) space of all finitely supported sequences $c:I\To \K$ with the inductive topology ($D$ is open in $\K^{(I)} \iff $ for every finite $J\subset I$ the set $D|_J=\{d\in \mathbb K^{J}: d=c|_J \text{ for some } c\in D\}$ is open in $\mathbb K^{J}$).

Consider the continuous bijection $U:X'\oplus  \K^{(I)}\To \overline{f(X')}$ given by $U(x',c)= f(x')+\sum_{i\in I}c(i)e_i$. We want to see that $U^{-1}:  \overline{f(X')}\To X'\oplus  \K^{(I)}$ is continuous. However $U^{-1}$ has the same graph as $U$, which is closed since $U$ is continuous.  The result we need is \cite[Theorem 14.7.1]{n-b}: {\em an operator from a Fr\'echet space to a webbed space whose graph is closed is continuous}. 

Banach (actually Fr\'echet) spaces are webbed \cite[Theorem 14.6.2]{n-b} and the product of two (even countably many) webbed spaces is again webbed. The space $\K^{(I)}$, being the strong dual of $\K^{I}$ (product topology), is also webbed (\cite[Theorem 14.6.4]{n-b}; this is the point where we use that $I$ is countable), and therefore 
$U^{-1}$ is continuous, which is possible only if $I$ is finite (otherwise $\K^{(I)}$ cannot be given a complete metric since it is the union of countably many closed sets with empty interior, see \cite[Section 11.8]{n-b}) and $f(X')$ is closed in $X$.
\end{proof}

The just proved result was suggested by \cite[Lemma~4.5]{cuchi}, which corresponds to the case where $\overline{f(X')}/f(X')$ is finite-dimensional.

We now tackle the problem of finding
  monomorphisms $f:X' \rightarrowtail X, g: Y'\rightarrowtail Y$ with dense range and separable codomain that are not isomorphic in $\W$.

To gauge the depths of our ignorance, consider the inclusions
$\imath:\ell_p\rightarrowtail \ell_q$
for  $1\leq p<q<\infty$: 
on the one hand we 
do not know whether there exist $(p,q)\neq (r,s)$ such that $\ell_p\rightarrowtail \ell_q$ and  $\ell_r\rightarrowtail \ell_s$ are not isomorphic in $\W$; on the other hand we have been unable to prove that any two of them are isomorphic!

It is relatively easy to see that some of these inclusions are not isomorphic in $\hM$ (take $\ell_1\rightarrowtail \ell_2$ and $\ell_2\rightarrowtail \ell_3$ if you want a really simple case) and that, in all cases,  $S(\imath)$ carries the trivial seminorm and has a Hamel basis whose cardinality is the continuum. By the way, can $\overline{f(X')}/f(X')$ consistently have a  Hamel basis whose cardinality is strictly smaller than the continuum?

\subsection*{Strictly singular objects}
An operator between Banach spaces is said to be strictly singular if it restriction to each infinite-dimensional closed subspace of its domain fails to be an isomorphism (linear homeomorphism onto its range), see \cite[Chapter 2, Section c]{l-t} for the basics. The definition applies in particular to the objects of $\W$. Assume 
$f:X' \rightarrowtail X$ is not strictly singular and take an infinite-dimensional closed subspace $A\subset X'$ witnessing it. Then we can ``cut-off'' the part of $f$ ``sitting on $A$'' since $f(A)$ is closed in $X$ and the right-square of the diagram
$$
\xymatrixcolsep{3pc}
\xymatrix{
A  \monicar[d]_{\text{restriction}} \ar[r]^-{\text{inclusion}} & X' \ar[r]^-{\text{quotient}} \monicar[d]^f & X'/A 
 \monicar[d]^-{\text{induced}}\\
f(A)  \ar[r]^{\text{inclusion}} & X \ar[r]^-{\text{quotient}} & X/f(A)\\
}
$$
is a pulation. This makes the hypotheses of the following statement less painful.

\begin{nada}\label{noisom} Let $K$ be an infinite-dimensional Banach space.
Assume $f:K \rightarrowtail \ell_1(I)$ is strictly singular and that $L$ is a Banach space such that every endomorphism of $K$ that factorizes through $L$ is strictly singular. Then no monomorphism
$g:L \rightarrowtail \ell_1(J)$ can be isomorphic to $f$ in $\W$.
\end{nada} 

The result applies when $K,L$ are different spaces in the family $\ell_p$ for $1\leq p\leq\infty$ or $c_0$. Note that each of these spaces has plenty of strictly  singular monomorphisms into $\ell_1$ with dense range: think of diagonal operators implemented by sequences converging to zero (and belonging to the dual space, of course).

\begin{proof}
Let $\phi: f\To g$ be an isomorphism in $\W$. By {\bf \ref{arisesasQ(u)}} there exist $u, v$ in $\hM$ such that $\phi=Q(u), \phi^{-1}=Q(v)$ and {\bf \ref{Qu=0}} implies that $vu={\bf I}_f$ in $\hM$: 
$$
\xymatrix{
K  \monicar[d]_f \ar[r]^{u'} & L \ar[r]^{v'} \monicar[d]^-g & K 
 \monicar[d]^-f\\
\ell_1(I) \ar[r]^{u} & \ell_1(J) \ar[r]^{v} & \ell_1(I)\\
}
$$
Let $r:\ell_1(I)\To K$ be an operator witnessing it, so that
$
vu={\bf I}_{\ell_1(I)}+fr $ and thus $ v'u'={\bf I}_{K}+rf.
$ 
Since $f$ is strictly  singular so is $rf$ and ${\bf I}_{K}+rf=v'u'$ cannot be strictly  singular, against the hypothesis. 
\end{proof}

During the proof we have obtained the following complement of {\bf \ref{arisesasQ(u)}}:

\begin{nada}
The objects $f:K \rightarrowtail \ell_1(I)$ and  $g:L \rightarrowtail \ell_1(J)$ are isomorphic in $\W$ if and only if they are isomorphic in $\hM$. 
 \end{nada}

\subsection*{A Lindenstrauss--Rosenthal theorem} The result alluded to in the heading states that if
$K, L$ are {\em closed} subspaces of $\ell_1$ whose corresponding quotients are isomorphic (and not isomorphic to $\ell_1$), then 
there is an automorphism of $\ell_1$ such that $L=U(K)$, see \cite[Theorem 2.f.8]{l-t}. 
In particular $U$ restricts to an isomorphism between $K$ and $L$. If we aim to obtain a version of this result in $\W$ we must realize that an essential ingredient of the proof of the Lindenstrauss--Rosenthal theorem is that each infinite-dimensional closed subspace of $\ell_1$, in particular $K$ and $L$, contains a subspace isomorphic to $\ell_1$ and complemented in $\ell_1$. Which does not makes a lot of sense in $\W$, does it? (If $f:K \rightarrowtail \ell_1$ has dense range, then $\Hom_{\W}(f,b)=0$ for all Banach spaces $b$, in particular if $b=\ell_1$.) We have the following  version of Lindenstrauss-Rosenthal's theorem in $\W$:

\begin{nada}\label{lr}
The objects $f:K \rightarrowtail \ell_1(I)$ and  $g:L \rightarrowtail \ell_1(J)$ are isomorphic in $\W$  if and only if there is an automorphism $U$ and an isomorphism $U'$ (both relative to the category of Banach spaces) forming a commutative square
$$
\xymatrixcolsep{3pc}
\xymatrix{
K \oplus\ell_1(J) \ar[r]^{U'} \ar[d]_{f\times {\bf I}} &  \ell_1(I)\oplus L \ar[d]^{{\bf I} \times g}\\  
\ell_1(I) \oplus \ell_1(J) \ar[r]^{U } &  \ell_1(I)\oplus\ell_1(J)
}
$$ 
\end{nada}

We are aware that this looks a bit tautological. 
We challenge the reader to write down a complete proof.

In particular $f$ and $g$ cannot be isomorphic in $\W$ if $K$ has the AP / BAP / DDP / RNP / Schur property / is weakly sequentially complete / is complemented in its bidual / has cotype $q\in[2,\infty)$ / (fill in with any property of $\ell_1(I)$ that passes to complemented subspaces) / and $L$ does not.

We are aware that the above welter of acronyms may be unintelligible for readers unfamiliar with (intermediate) Banach space theory. Instead of trying to explain these definitions here, which would be useless, we refer the interested reader to Albiac--Kalton's book where they can be studied with minimal prerequisites:
\begin{itemize}
\item AP stands for Grothendieck's ``approximation property'' and BAP for the ``bounded approximation property'' see \cite[Definitions 1.4--1.5]{AK}.

\item DPP is the Dunford--Pettis property \cite[Definition 5.4.3]{AK}, also isolated by 
Grothendieck.

\item RNP is the Radon--Nikod\'ym property, see \cite[Definition 5.5.1]{AK}.

\item A Banach space has the Schur property if weakly convergent {\em sequences} converge in norm, see \cite[Definition 2.3.4]{AK}. The meaning of ``weakly sequentially complete'' and ``complemented in its bidual'' are self-explanatory, see \cite[Definition 2.3.4 and p. 45]{AK} just in case.

\item The notion of cotype is an invention of Hoffmann-J\o rgensen, see \cite[Definition 6.2.10 (b)]{AK}.

\end{itemize}

Admittedly, {\bf \ref{noisom}} and {\bf \ref{lr}} are quite unsatisfactory because they depend on the domains of $f$ and $g$ being very different. It would be delightful to have counterexamples with the same domain. The final section of the paper proposes a convenient setting to address this issue.

\subsection*{Enriched $\Ext ^n_{\bf B}$}
Quoting Mac Lane \cite[Chapter XII, first paragraph of \S9]{mac}, ``A standard method [to derive functors] is: Take a resolution, apply a covariant functor
$T$, take the homology of the resulting complex. This gives a
connected sequence of functors, called the derived functors of $T$.''

In the end, the functors $\Ext^n_{\bf B}$ are the derived functors of $\Hom$ in the most orthodox sense of the word in $\W$, that is, when $\Hom$ is treated as a functor from $\W$ to $\bf V$, derived as such and one restricts the result to $\B$, %
see  \cite[Part II, Section 21]{HHI}. 


If $a,b$ are Banach spaces, then so is $\Hom_{\B}(a,b)$ and, therefore, the ``object'' $\Ext_{\B}^n(a,b)$ appears as the cokernel space of
$
\varkappa^*: \Hom_{\B}(p_n,b)\To \Hom_{\B}(k,b)
$
either in $\bf V$ or in $\bf sB$ (see \cite[Section 4.5]{hmbst}).
One can consider an ``enriched version'' of  $\Ext_{\B}^n$ taking the cokernel morphism (which is an object of $\W$). The resulting functor takes values in $\W$ and the traditional $\bf sB$-valued version appears after composing with $S$. Related ideas, under a different perspective and in a different context, appear in the first section of the introductory chapter of \cite{buhler}, see especially the theorem on p.~x and the corollary on p.~xii.

\subsection*{Quotient Hilbert / quasi-Banach spaces}


Let $\bf H$ be the category of Hilbert spaces and their operators: this is a full subcategory of $\B$ from which it inherits most universal constructions (kernels and cokernels, pullbacks and pushouts), nice properties (the quasi abelian character) and shortcomings (it fails to be abelian and does not have injective or projective objects, apart from zero).

Let ${\bf MH}$ be the category of monics in $\bf H$ and ${\bf hMH}$ the homotopic version. Both   ${\bf MH}$  and ${\bf hMH}$
contain $\bf H$ as a full subcategory.

The main structural difference between  ${\bf hMH}$ and its Banach relative $\hM$ stems from the fact that all short exact sequences of Hilbert spaces split. Since {\bf \ref{pul=isom}} holds in ${\bf hMH}$ we see that pulations are automatically isomorphisms in ${\bf hMH}$ and so ${\bf hMH}$ is abelian --- no further localization is required!

It follows that  $f:H' \rightarrowtail H$ is isomorphic to a genuine Hilbert space if and  only if $f(H')$ is closed in $H$ if and  only if $f$ is projective in ${\bf hMH}$. There are no injective objects in ${\bf hMH}$, apart from zero (see {\bf \ref{noinj}}, although in this case the proof is much easier).

Moreover, if $f$ lives in ${\bf hMH}$ the short exact sequence (\ref{SEX})
splits and so $f$ can be written as the direct sum of a Hilbert space (the ``quotient'' part) and a monomorphism with dense range (the ``subspace'' part) to which {\bf \ref{uncountable}} applies.

The content of this paper flows naturally into the following problems in classical operator theory:

\medskip

\noindent$\bigstar$ 
 Under which conditions are two injective endomorphisms of $\ell_1$ (respectively, of $\ell_2$) with dense range  isomorphic in $\hM$   (respectively, in ${\bf hMH}$)?

\medskip 

It is quite possible, however, that the following question is much more thrilling:

\medskip

\noindent$\bigstar$ 
What can be said about the category of ``quotient quasi-Banach spaces''?

\end{document}